\documentclass[12pt]{amsart}
\usepackage[margin=1.3in]{geometry}
\usepackage[all,cmtip]{xy}

\theoremstyle{plain}
\newtheorem{theorem}{Theorem}
\newtheorem{prop}[theorem]{Proposition}
\newtheorem{corollary}[theorem]{Corollary}
\newtheorem{lemma}[theorem]{Lemma}

\theoremstyle{definition}
\newtheorem{defn}[theorem]{Definition}
\newtheorem{rmk}[theorem]{Remark}

\numberwithin{theorem}{section}

\newcommand{\spec}{\operatorname{Spec}}
\newcommand{\ord}{\operatorname{ord}}
\newcommand{\codim}{\operatorname{codim}}
\newcommand{\Cont}{\operatorname{Cont}}

\begin{document}

\title{Jet schemes of determinantal varieties}

\author[C. Yuen]{Cornelia Yuen} 
\address{SUNY Potsdam, Department of Mathematics, 44 Pierrepont
Avenue, Potsdam, NY 13676, USA} 
\email{{\tt yuenco@potsdam.edu}}

\begin{abstract}
    This article studies the scheme structure of the jet schemes of
    determinantal varieties.  We show that in general, these jet
    schemes are not irreducible.  In the case of the determinantal
    variety $X$ of $r \times s$ matrices of rank at most one, we give
    a formula for the dimension of each of the components of its jet
    schemes.  As an application, we compute the log canonical
    threshold of the pair $(\mathbb{A}^{rs},X)$.
\end{abstract}

\thanks{This is part of the author's PhD thesis at the University of
Michigan under the supervision of Karen Smith.  The author was
supported by her advisor's NSF grant $0070722$ and the department's
NSF RTG grant $0502170$.} 

\thanks{The author would like to thank Mircea Musta\c{t}\v{a} for
suggesting this problem, and for many valuable discussions.}

\maketitle

\section{Introduction}

Let $X$ be a scheme of finite type over an algebraically closed field
$k$ of characteristic zero.  An arc on an algebraic variety $X$ is an
``infinitesimal curve'' on it.  Formally, this is a morphism, defined
over $k$, from the curve germ scheme $\spec k[[t]]$ into $X$.  The set
of all arcs of $X$ carries the structure of a scheme, called the arc
space of $X$ and denoted by $\mathcal{J}_{\infty}(X)$.

An $m$-jet on $X$ is a truncated arc on $X$, that is, a $k$-morphism
$$\spec k[t]/(t^{m+1}) \rightarrow X.$$ The set of all $m$-jets on $X$
also forms a scheme in a natural way.  This is the $m^{th}$ jet scheme
$\mathcal{J}_{m}(X)$.

The surjection $k[t]/(t^{m+1}) \rightarrow k[t]/(t^{m})$ induces a
morphism $\pi^{m}_{m-1}:\mathcal{J}_{m}(X) \rightarrow
\mathcal{J}_{m-1}(X)$, and composition gives a morphism
$\pi_{m}:\mathcal{J}_{m}(X) \rightarrow \mathcal{J}_{0}(X)=X$. Taking 
the inverse limit of these jet schemes gives the arc space 
$\mathcal{J}_{\infty}(X)$.

In \cite{MustataLCI}, Musta\c{t}\v{a} proved that a locally complete
intersection variety has rational singularities if and only if all its
jet schemes are irreducible.  Combining this result with the fact that
determinantal varieties have rational singularities, we know that a
determinantal variety of singular square matrices is a hypersurface
and so all its jet schemes are irreducible.  Though determinantal
varieties always have rational singularities, they are rarely complete
intersections, so Musta\c{t}\v{a}'s theorem is not applicable to more
general determinantal varieties.  This leads to the natural question:
Are the jet schemes of all determinantal varieties irreducible?

In this paper, we will see that the answer to this question is no.
In particular, we will show that the second jet scheme and all the odd
jet schemes of essentially all determinantal varieties are reducible.
In the special case of the variety of matrices of rank at most one, we
will give a formula for the number of irreducible components of
its jet schemes and their dimensions. 

\section{Background and notation}

Let $X_{c} \subseteq \mathbb{A}^{rs}$ be the determinantal variety of $r
\times s$ matrices of rank at most $c$, that is, $X_{c}$ is defined by
$I_{c+1}$, the ideal of $(c+1)$-minors\footnote{We define a $d$-minor
of a matrix to be a $d \times d$ subdeterminant of that matrix.} of a
generic $r \times s$ matrix $(x_{ij})$.  So an $m$-jet of $X_{c}$
corresponds to a $k$-algebra homomorphism
\begin{align*}
    \phi: k[x_{ij}] &\rightarrow k[t]/(t^{m+1})\\
    x_{ij} &\mapsto
    x_{ij}^{(0)}+x_{ij}^{(1)}t+\ldots+x_{ij}^{(m)}t^{m},
\end{align*}
subject to the condition $\phi(I_{c+1})=0$.

For a generator $\triangle \in I_{c+1}$, we write
$$\phi(\triangle)=\triangle^{(0)}+\triangle^{(1)}t+\ldots+\triangle^{(m)}t^{m}$$
where $\triangle^{(l)} \in k[x_{ij}^{(k)}]$.  These polynomials
$\triangle^{(l)}$ can be obtained by calculating the corresponding
$(c+1)$-minors of the matrix whose entries are
$x_{ij}^{(0)}+x_{ij}^{(1)}t+\ldots+x_{ij}^{(m)}t^{m}$, and extracting
the coefficients of each power of $t$.  Therefore,
$\mathcal{J}_{m}(X_{c})$ is the closed subscheme of
$\mathcal{J}_{m}(\mathbb{A}^{rs})=\spec k[x_{ij}^{(k)}]$ defined by
the ideal $(\triangle^{(0)},\triangle^{(1)},\ldots,\triangle^{(m)}:
\triangle \text{'s are }(c+1)\text{-minors of }(x_{ij}))$.

Before we present the results, let us outline the general strategy.

A general irreducible scheme $X$ over a field $k$ is the disjoint
union of its singular locus $X^{sing}$ and its smooth locus
$X^{reg}$. So to understand the jet scheme $\mathcal{J}_{m}(X)$, we
can study the preimage of these two loci under the natural
projection $\pi_{m}:\mathcal{J}_{m}(X) \rightarrow X$. Since
$X^{reg}$ is smooth and irreducible, $\pi_{m}$ is an affine bundle
over $X^{reg}$, and therefore $\overline{\pi_{m}^{-1}(X^{reg})}$ is
an irreducible component of $\mathcal{J}_{m}(X)$ of dimension
$(m+1)\dim X$. On the other hand, $\pi_{m}^{-1}(X^{sing})$ is a
closed subset of $\mathcal{J}_{m}(X)$. So it will contribute
components to $\mathcal{J}_{m}(X)$ except when it is contained in
$\overline{\pi_{m}^{-1}(X^{reg})}$ (as is the case when $X$ is a
local complete intersection with rational singularities).
Therefore, if we want to show $\mathcal{J}_{m}(X)$ is reducible, it
suffices to show $\dim \pi_{m}^{-1}(X^{sing}) \geq \dim
\overline{\pi_{m}^{-1}(X^{reg})}$.

\begin{rmk}\label{Converse}
    There are examples when $\mathcal{J}_{m}(X)$ is reducible even
    though $\dim \pi_{m}^{-1}(X^{sing})< \dim
    \overline{\pi_{m}^{-1}(X^{reg})}$, see Remark
    \ref{ConverseExample}.
\end{rmk}

\section{Odd jet schemes are reducible}

Our first result says that all the odd jet schemes of essentially all 
determinantal varieties are reducible. More precisely,

\begin{theorem}\label{OddJets}
    Let $X_{c}$ be the determinantal variety of $r \times s$ matrices
    of rank at most $c$, $c \geq 1$.  If $r,s \geq c+2$ and $r+s \geq
    2c+5$, then the jet scheme $\mathcal{J}_{m}(X_{c})$ is reducible 
    for $m$ odd.
\end{theorem}

To prove this theorem, we need the following general result:

\begin{lemma}\label{Dimension}
    Let $X$ be a smooth scheme and $Z \subseteq Y \subseteq X$ so that
    $Y$ is closed in $X$ and $Z$ is closed in $Y$.  Then $\dim
    \pi_{mp-1}^{-1}(Z) \geq p \dim \pi_{m-1}^{-1}(Z)$ for all $m,p
    \geq 1$, where $\pi_{i}:\mathcal{J}_{i}(Y) \rightarrow Y$ are the
    natural projections as described in the Introduction.
\end{lemma}

To prove this lemma, we will need a general fact (Proposition
\ref{ContactLoci}, below) about contact loci whose statement and proof
were provided by Musta\c{t}\v{a}.  But first let us look at some
basics of contact loci (see for example \cite[\S 2.4 and \S
5]{Blickle}).

\begin{defn}
    Let $Y \subseteq X$ be a closed subscheme of a smooth scheme $X$ 
    defined by the sheaf of ideals $\mathcal{I}_{Y}$. We define the 
    function $$\ord_{Y}:\mathcal{J}_{\infty}(X) \rightarrow 
    \mathbb{N} \cup \{\infty\}$$ sending an arc 
    $\gamma:\mathcal{O}_{X} \rightarrow k[[t]]$ of $X$ to the order 
    of vanishing of $\gamma$ along $Y$, that is, $\ord_{Y}(\gamma)$ 
    is the integer $e$ such that the ideal $\gamma(\mathcal{I}_{Y}) 
    \subseteq k[[t]]$ is exactly the ideal $(t^{e})$.
\end{defn}

\begin{defn}
    With $Y \subseteq X$ as above, the contact locus $\Cont^{\geq
    m}(Y)$ is the set of arcs of $X$ whose order of vanishing along
    $Y$ is at least $m$.  In other words, $$\Cont^{\geq m}(Y)=\{\gamma
    \in \mathcal{J}_{\infty}(X)|\ord_{Y}(\gamma) \geq
    m\}=\mu_{m-1}^{-1}(\mathcal{J}_{m-1}(Y))$$ where
    $\mu_{m-1}:\mathcal{J}_{\infty}(X) \rightarrow
    \mathcal{J}_{m-1}(X)$ is the natural projection.  In this case, 
    we define 
    \begin{align*}
	\codim(\mathcal{J}_{\infty}(X),\Cont^{\geq m}(Y)) &:= \codim
	(\mathcal{J}_{m-1}(X),\mu_{m-1}(\Cont^{\geq m}(Y)))\\
	&= \codim(\mathcal{J}_{m-1}(X),\mathcal{J}_{m-1}(Y)).
    \end{align*}
    We can also define
    \begin{align*}
	\Cont^{m}(Y) &= \Cont^{\geq m}(Y) \setminus \Cont^{\geq
	m+1}(Y)\\
	&= \mu_{m-1}^{-1}(\mathcal{J}_{m-1}(Y)) \setminus
	\mu_{m}^{-1}(\mathcal{J}_{m}(Y)),
    \end{align*}
    and
    $\codim(\mathcal{J}_{\infty}(X),\Cont^{m}(Y))=
    \codim(\mathcal{J}_{\infty}(X),\Cont^{\geq m}(Y))$, because 
    $\Cont^{m}(Y)$ is a non-empty open subset of $\Cont^{\geq m}(Y)$.
\end{defn}

To understand the statement of Proposition \ref{ContactLoci}, we 
also need to recall the definition of a log resolution.  

\begin{defn}\label{LogResolution}
    Let $X$ be a smooth variety over a field $k$ of characteristic 
    zero, and $Y \subseteq X$ a closed subscheme. A map 
    $f:\widetilde{X} \rightarrow X$ is called a log resolution of the 
    pair $(X,Y)$ if $f$ is proper birational such that
    \begin{enumerate}
	\item $\widetilde{X}$ is smooth and
	\item writing $F:=f^{-1}(Y)=\sum_{i=1}^{r} a_{i}D_{i}$ and 
	$K_{\widetilde{X}/X}=\sum_{i=1}^{r} k_{i}D_{i}$ for some 
	$a_{i}$, $k_{i} \in \mathbb{Q}$ and prime divisors $D_{i}$, 
	the divisors $F$, $K_{\widetilde{X}/X}$ and 
	$F+K_{\widetilde{X}/X}$ have simple normal crossing support.
    \end{enumerate}	
\end{defn}

\begin{rmk}
    The existence of log resolutions follows from Hironaka's 
    resolution of singularities \cite{Hironaka}. Note also that one 
    can find a simultaneous log resolution for any finite set 
    $Y_{1},\ldots,Y_{r}$ of closed subschemes of $X$.
\end{rmk}

\begin{prop}\label{ContactLoci}
    Let $\{Y_{i}\}_{i=1}^{r}$ be a collection of closed subschemes of
    a smooth scheme $X$.  Let $f:X' \rightarrow X$ be a log resolution
    of $Y_{1},\ldots,Y_{r}$, and denote $$f^{-1}(Y_{i})=\sum_{j=1}^{n}
    a_{ij}E_{j} \hspace{0.2cm}\text{and}\hspace{0.2cm}
    K_{X'/X}=\sum_{j=1}^{n} k_{j}E_{j},$$ where $\sum E_{j}$ is a
    simple normal crossing divisor.  Then $$\codim \left(
    \bigcap_{i=1}^{r} \Cont^{\geq
    m_{i}}(Y_{i})\right)=\min_{\underline{\nu} \in \mathbb{N}^{n}}
    \left\{\sum_{j=1}^{n} \nu_{j}(k_{j}+1)\left|\bigcap_{\nu_{j} \geq
    1} E_{j} \neq \emptyset \text{ \rm{and} }\sum_{j=1}^{n}
    a_{ij}\nu_{j} \geq m_{i} \forall i \right.\right\}.$$
\end{prop}

\begin{proof}
    Let $C=\bigcap_{i=1}^{r} \Cont^{\geq m_{i}}(Y_{i})$ and
    $f_{\infty}:\mathcal{J}_{\infty}(X') \rightarrow
    \mathcal{J}_{\infty}(X)$ be the induced map on arc spaces.  Then
    \begin{align*}
	f_{\infty}^{-1}(C) &= \{\gamma \in
	\mathcal{J}_{\infty}(X')|f_{\infty}(\gamma) \in C\}\\
	&= \{\gamma \in
	\mathcal{J}_{\infty}(X')|\ord_{Y_{i}}(f_{\infty}(\gamma)) \geq
	m_{i} \text{ for }1 \leq i \leq r\}\\
	&= \{\gamma \in \mathcal{J}_{\infty}(X')|\sum_{j}
	a_{ij}\ord_{E_{j}}(\gamma) \geq m_{i} \text{ for }1 \leq i
	\leq r\}.
    \end{align*}
    Now for an arc $\gamma$ of $X'$ with order of vanishing $\nu_{j}$
    along $E_{j}$ for all $j$, $\gamma \in f_{\infty}^{-1}(C)$ if and
    only if $\sum_{j=1}^{n}a_{ij}\nu_{j} \geq m_{i}$ for $1 \leq i
    \leq r$.  So we have
    $$f_{\infty}^{-1}(C)=\bigsqcup_{\underline{\nu}}
    \left(\bigcap_{j=1}^{n} \Cont^{\nu_{j}}(E_{j})\right)$$ where the
    disjoint union is taken over all $\underline{\nu} \in
    \mathbb{N}^{n}$ such that $\sum_{j} a_{ij}\nu_{j} \geq m_{i}$ for
    all $i=1,\ldots,r$.  This implies $$C=\bigsqcup_{\underline{\nu}}
    f_{\infty} \left(\bigcap_{j=1}^{n}
    \Cont^{\nu_{j}}(E_{j})\right),$$ and thus $$\codim
    C=\min_{\underline{\nu}} \left\{\codim
    f_{\infty}\left(\bigcap_{j=1}^{n}
    \Cont^{\nu_{j}}(E_{j})\right)\right\}$$
    where the minimum is taken over all $\underline{\nu} \in
    \mathbb{N}^{n}$ such that $\bigcap_{\nu_{j} \geq 1}E_{j} \neq
    \emptyset$ and $\sum_{j} a_{ij}\nu_{j} \geq m_{i}$ for all
    $i=1,\ldots,r$.  By \cite[Theorem 2.1]{EinLazarsfeldMustata},
    $$\codim
    f_{\infty}\left(\bigcap_{j=1}^{n}\Cont^{\nu_{j}}(E_{j})\right)=
    \sum_{j=1}^{n}\nu_{j}(k_{j}+1).$$ As a result, $$\codim
    C=\min_{\underline{\nu} \in \mathbb{N}^{n}} \left\{\sum_{j=1}^{n}
    \nu_{j}(k_{j}+1)\left|\bigcap_{\nu_{j} \geq 1} E_{j} \neq
    \emptyset \text{ and }\sum_{j=1}^{n} a_{ij}\nu_{j} \geq m_{i}
    \forall i \right.\right\}.$$
\end{proof}
    
\begin{proof}[Proof of Lemma \ref{Dimension}]
    Taking $Y_{1}=Y$, $Y_{2}=Z$, $m_{1}=m$ and $m_{2}=1$ in
    Proposition \ref{ContactLoci}, we have
    \begin{align*}
	&\codim \left(\Cont^{\geq m}(Y)\cap \mu^{-1}(Z)\right)\\
	= &\min_{\underline{\nu} \in \mathbb{N}^{n}}
	\left\{\sum_{j=1}^{n} \nu_{j}(k_{j}+1)\left|\bigcap_{\nu_{j}
	\geq 1} E_{j} \neq \emptyset \text{ and }\sum_{j=1}^{n}
	a_{j}\nu_{j} \geq m \right.\right\},
    \end{align*}
    where $\mu:\mathcal{J}_{\infty}(Y) \rightarrow Y$ is the natural
    projection.  Say $\underline{\nu} \in \mathbb{N}^{n}$ is the
    $n$-tuple achieving the minimum value above.  Then the $n$-tuple
    $p \underline{\nu}$ satisfies the conditions $\bigcap_{p\nu_{j}
    \geq 1}E_{j} \neq \emptyset$ and $\sum_{j=1}^{n}a_{j}(p\nu_{j})
    \geq mp$.  So we obtain
    \begin{equation}
	\codim \left(\Cont^{\geq mp}(Y)\cap \mu^{-1}(Z)\right) \leq p
	\cdot \codim \left(\Cont^{\geq m}(Y)\cap \pi^{-1}(Z)\right).
	\tag{$\dagger$}
    \end{equation}
    On the other hand,
    \begin{align*}
	\codim \left(\Cont^{\geq m}(Y)\cap \pi^{-1}(Z)\right) &=
	\codim(\mathcal{J}_{m-1}(X),\mathcal{J}_{m-1}(Y) \cap
	\pi_{m-1}^{-1}(Z))\\
	&= \dim \mathcal{J}_{m-1}(X)-\dim \pi_{m-1}^{-1}(Z)\\
	&= m\dim X-\dim \pi_{m-1}^{-1}(Z).
    \end{align*}
    Therefore, the inequality ($\dagger$) gives us $$mp\dim X-\dim
    \pi_{mp-1}^{-1}(Z) \leq p(m\dim X-\dim \pi_{m-1}^{-1}(Z)),$$ or
    equivalently, $\dim \pi_{mp-1}^{-1}(Z) \geq p \dim
    \pi_{m-1}^{-1}(Z)$, as desired.
\end{proof}

\begin{proof}[Proof of Theorem \ref{OddJets}]
    Because $m$ is odd, we can write $m=2p-1$ for some $p \geq 1$.
    According to our discussion prior to Remark \ref{Converse}, to
    show that $\mathcal{J}_{m}(X_{c})$ is reducible, it suffices to
    show $\dim \pi_{m}^{-1}(X_{c}^{sing}) \geq \dim
    \overline{\pi_{m}^{-1}(X_{c}^{reg})}$.  By Lemma \ref{Dimension},
    we have $\dim \pi_{m}^{-1}(X_{c}^{sing}) \geq p \dim
    \pi_{1}^{-1}(X_{c}^{sing})$.  So we need to compute the dimension
    of $\pi_{1}^{-1}(X_{c}^{sing})$.

    Recall that $X_{c}^{sing}$ is the set of all $r \times s$ matrices
    of rank at most $c-1$, and hence is the subvariety of $X_{c}$
    defined by the $c$-minors of the generic matrix $(x_{ij})$.  Let
    $B=(x_{ij}^{(0)})$ be an $r \times s$ matrix of indeterminates.
    Then
    \begin{align*}
	\pi_{1}^{-1}(X_{c}^{sing}) &= \spec
	k[x_{ij}^{(0)},x_{ij}^{(1)}]/(J_{1}(X_{c})+(c\text{-minors
	of }B))\\
	&= \spec k[x_{ij}^{(0)},x_{ij}^{(1)}]/(c\text{-minors of
	}B),
    \end{align*}
    since $J_{1}(X_{c}) \subseteq (c\text{-minors of }B)$.  So
    $\pi_{1}^{-1}(X_{c}^{sing})$ has dimension $rs+(c-1)(r+s-c+1)$.

    Since $X_{c}$ has dimension $c(r+s-c)$, the component of
    $\mathcal{J}_{m}(X_{c})$ over the smooth part of $X_{c}$ has dimension
    $(m+1)\dim X_{c}$, which is $2pc(r+s-c)$.  Using our hypotheses
    that $r,s \geq c+2$ and $r+s \geq 2c+5$, we see that
    $(r-c-1)(s-c-1) \geq 2$, or equivalently, $$rs+(c-1)(r+s-c+1) \geq
    2c(r+s-c).$$ Therefore, $\dim \pi_{m}^{-1}(X_{c}^{sing}) \geq \dim
    \overline{\pi_{m}^{-1}(X_{c}^{reg})}$.  Thus, the preimage of the 
    singular locus of $X_{c}$ in $\mathcal{J}_{m}(X_{c})$ has 
    dimension at least as large as the dimension of 
    $\mathcal{J}_{m}(X_{c})$ over the generic point of $X_{c}$. As 
    explained in the paragraph prior to Remark \ref{Converse}, 
    therefore, $\mathcal{J}_{m}(X_{c})$ is not irreducible.
\end{proof}

\begin{rmk}
    The argument above tells us that $$\dim
    \mathcal{J}_{1}(X_{c})=rs+(c-1)(r+s-c+1)$$ provided that $c \geq
    1$, $r,s \geq c+2$ and $r+s \geq 2c+5$.
\end{rmk}

\begin{rmk}\label{ConverseExample}
    Although our dimension analysis draws no conclusion to the
    irreducibility of odd jet schemes of $X_{c}$ when $c \geq 1$ and
    $r=s=c+2$, there are examples demonstrating that
    $\mathcal{J}_{m}(X_{c})$ is still reducible.  For example, take
    $X$ to be the variety of $3 \times 3$ matrices of rank at most
    one.  Then the proof of Theorem \ref{OddJets} tells us that $$\dim
    \pi_{1}^{-1}(X^{sing})=9<10=\dim
    \overline{\pi_{1}^{-1}(X^{reg})}.$$ However, a Macaulay
    calculation says the ideal $J_{1}(X) \subset
    k[x_{ij}^{(0)},x_{ij}^{(1)}]$ for $1 \leq i,j \leq 3$, has two
    minimal primes: $$J_{1}(X)+(x_{11}^{(1)}x_{22}^{(1)}x_{33}^{(1)}
    -x_{11}^{(1)}x_{32}^{(1)}x_{23}^{(1)}
    -x_{12}^{(1)}x_{21}^{(1)}x_{33}^{(1)}
    +x_{12}^{(1)}x_{31}^{(1)}x_{23}^{(1)}
    +x_{13}^{(1)}x_{21}^{(1)}x_{32}^{(1)}
    -x_{13}^{(1)}x_{31}^{(1)}x_{22}^{(1)})$$ and
    $$(x_{11}^{(0)},x_{12}^{(0)},x_{13}^{(0)},x_{21}^{(0)},
    x_{22}^{(0)},x_{23}^{(0)},x_{31}^{(0)},x_{32}^{(0)},x_{33}^{(0)}).$$
    That is, the jet scheme $\mathcal{J}_{1}(X)$ has two irreducible
    components and is therefore reducible.
\end{rmk}

\section{Second jet scheme is reducible}

In this section, we investigate the second jet scheme of a
determinantal variety.

We denote by $A_{k}$ the $r \times s$ matrix with a $k \times k$
identity submatrix in the upper left corner and zero entries
everywhere else.

\begin{theorem}\label{SecondJets}
    Let $X_{c}$ be the variety of $r \times s$ matrices of rank at
    most $c$, $c \geq 2$.  Suppose $r,s \geq c+2$ and $r+s \geq 2c+6$.
    Then $\mathcal{J}_{2}(X_{c})$ is reducible.
\end{theorem}
\begin{proof}
    As outlined in the paragraph prior to Remark \ref{Converse}, our
    goal is to show that $\dim \pi_{2}^{-1}(X_{c}^{sing}) \geq \dim
    \overline{\pi_{2}^{-1}(X_{c}^{reg})}$.  First notice that the
    singular locus of $X_{c}$ can be stratified, according to the rank
    of the singular points.  So to understand
    $\pi_{2}^{-1}(X_{c}^{sing})$, we can study the preimage of
    matrices of a fixed rank.  Notice that the group $GL(r) \times
    GL(s)$ acts transitively on matrices of a fixed rank.  So when we
    consider the fiber of $\pi_{2}$ over a singular point of rank
    $k<c$, we may pick the representative $A_{k}$.  Now we describe
    the fiber over $A_{k}$.

    The ideal $J_{2}(X_{c})$ is homogeneous.  Its generators are the
    coefficients of $t^{0}$, $t^{1}$ and $t^{2}$ in the $(c+1)$-minors
    of the $r \times s$ matrix $$\left(
    x_{ij}^{(0)}+x_{ij}^{(1)}t+x_{ij}^{(2)}t^{2} \right)_{\substack{1
    \leq i \leq r\\ 1 \leq j \leq s}}.$$
    Thus, every generator has degree $c+1$ and every term of each of
    its generators has at least $c-1$ variables of the form
    $x_{ij}^{(0)}$'s.  So over the singular point $A_{k}$, for $k \leq
    c-2$, we have $$\pi_{2}^{-1}(A_{k})=\mathcal{J}_{2}(X_{c}) \times
    \spec k[x_{ij}^{(0)}]/m_{A_{k}}$$ where $m_{A_{k}}$ is the maximal
    ideal of the point $A_{k} \in X_{c-2}$.  Since $$m_{A_{k}}=\left(
    x_{pp}^{(0)}-1:1 \leq p \leq k, x_{mn}^{(0)}:m=n>k \text{ or } m
    \neq n\right),$$ we have
    \begin{align*}
	&\hspace{0.3in} k[x_{ij}^{(0)},x_{ij}^{(1)},x_{ij}^{(2)}]/
	\left( J_{2}(X_{c})+\left( x_{pp}^{(0)}-1:1 \leq p \leq k,
	x_{mn}^{(0)}:m=n>k \text{ or } m \neq n\right)\right)\\
	&\cong k[x_{ij}^{(0)},x_{ij}^{(1)},x_{ij}^{(2)}]/ \left(
	x_{pp}^{(0)}-1:1 \leq p \leq k, x_{mn}^{(0)}:m=n>k \text{ or }
	m \neq n\right).
    \end{align*}

    This means that the fiber over any singular point of rank at most
    $c-2$ is isomorphic to $\mathbb{A}^{2rs}$, and therefore $$\dim
    \pi_{2}^{-1}(X_{c-2})=2rs+(c-2)(r+s-c+2).$$

    Over the singular point $A_{c-1}$, a surviving term of a generator
    of $J_{2}(X_{c})$ has the form $$x_{11}^{(0)}x_{22}^{(0)}\cdots
    x_{c-1,c-1}^{(0)}x_{ij}^{(1)}x_{kl}^{(1)}$$ where $x_{ij}^{(1)}$
    and $x_{kl}^{(1)}$ are two distinct entries in the lower right
    $(r-c+1) \times (s-c+1)$ submatrix of the $r \times s$ matrix
    $(x_{pq}^{(1)})$.  So over $A_{c-1}$, we obtain
    \begin{align*}
	&\hspace{0.3in} k[x_{ij}^{(0)},x_{ij}^{(1)},x_{ij}^{(2)}]/
	\left(J_{2}(X_{c})+\left( x_{pp}^{(0)}-1:1 \leq p \leq c-1,
	x_{mn}^{(0)}:m=n \geq c \text{ or } m \neq n\right)\right)\\
	&\cong k[x_{ij}^{(1)},x_{ij}^{(2)}]/(2\text{-minors of
	the matrix } (x_{kl}^{(1)})_{\substack{1 \leq k \leq r-c+1\\ 1 \leq l
	\leq s-c+1}}).
    \end{align*}
    This implies that the fiber over $A_{c-1}$ has dimension
    $$(r+s-2c+1)+(rs-(r-c+1)(s-c+1))+rs.$$  As a result, the preimage 
    of the set of rank $c-1$ matrices under the map $\pi_{2}$ has 
    dimension $$r+s-2c+1+2rs-(r-c+1)(s-c+1)+(c-1)(r+s-c+1).$$

    Now, to compare $$\dim \pi_{2}^{-1}(X_{c}^{sing})=\max \{\dim
    \pi_{2}^{-1} (X_{c-2}), \dim \pi_{2}^{-1}(X_{c-1} \setminus
    X_{c-2})\}$$ and $$\dim
    \overline{\pi_{2}^{-1}(X_{c}^{reg})}=3c(r+s-c),$$ we observe that
    $$\dim \pi_{2}^{-1}(X_{c-2}) \geq \dim
    \overline{\pi_{2}^{-1}(X_{c}^{reg})} \text{ if and only if }
    (r-c-1)(s-c-1) \geq 3$$ and $$\dim \pi_{2}^{-1}(X_{c-1} \setminus
    X_{c-2}) \geq \dim \overline{\pi_{2}^{-1}(X_{c}^{reg})} \text{ if
    and only if } (r-c-1)(s-c-1) \geq 2.$$ But our hypotheses $r,s
    \geq c+2$ and $r+s \geq 2c+6$ are equivalent to the condition
    $(r-c-1)(s-c-1) \geq 3$.  This completes the proof that
    $\mathcal{J}_{2}(X_{c})$ is reducible with our assumptions
    on $r$, $s$ and $c$.
\end{proof}

\section{Varieties of matrices of rank at most one}

While the analysis on the scheme structure of the jet schemes of a
general determinantal variety remains incomplete, the case when $t=1$
is much better understood.  This is largely due to the fact that the
singular locus of this type of determinantal varieties is an isolated 
origin and that we have a very nice description of the preimage
of this singular set under the map $\pi_{m}$.

Musta\c{t}\v{a} showed that the higher jet schemes of the
determinantal variety of $2 \times n$ matrices of rank at most one are
all irreducible \cite[Example 4.7]{MustataLCI}.  However, this result
does not hold for larger matrices.  In fact, we have a complete 
understanding of the number of components of the jet schemes in this 
case, and a formula for the dimension of each of the components.

\begin{theorem}\label{2x2minors}
    Let $X$ be the variety of $r \times s$ matrices of rank at most
    one.  Assume $r>s \geq 3$.  Then
    $\mathcal{J}_{m}(X)$ has precisely $\lfloor \frac{m+1}{2} \rfloor
    +1$ irreducible components and these components have dimensions
    $qrs+(m+1-2q)\dim X$ where $q=0,\ldots,\lfloor \frac{m+1}{2}
    \rfloor$.  In particular, the dimension of $\mathcal{J}_{m}(X)$ is
    $=\lfloor \frac{m+1}{2} \rfloor rs+(m \text{ mod }2)\dim X$.
\end{theorem}

As preparation for the proof of this theorem, let us first examine the
preimage of the origin under the natural projection
$\pi_{m}:\mathcal{J}_{m}(X) \rightarrow X$ in a slightly more 
general context.

\begin{prop}\label{PreimageOfOrigin}
    Let $X \subseteq \mathbb{A}^{n}$ be a closed subscheme defined
    over $k$ by a set of homogeneous polynomials, all of the same
    degree $d$.  Let $\pi_{m}:\mathcal{J}_{m}(X) \rightarrow X
    \subseteq \mathbb{A}^{n}$ be the natural surjection.  Then for $m
    \geq d \geq 2$, $$\pi_{m}^{-1}(0) \cong \mathcal{J}_{m-d}(X)
    \times \mathbb{A}^{n(d-1)}.$$
\end{prop}
\begin{proof}
    Let $I$ be the defining ideal of $X$.  Then an $m$-jet of $X$
    corresponds to a ring homomorphism
    \begin{align*}
	k[x_{1},\ldots,x_{n}]/I &\rightarrow k[t]/(t^{m+1})\\
	x_{i} &\mapsto
	x_{i}^{(0)}+x_{i}^{(1)}t+\ldots+x_{i}^{(m)}t^{m}
    \end{align*}
    where $x_{i}^{(j)} \in k$ are arbitrary.  This $m$-jet lies in
    $\pi_{m}^{-1}(0)$ if and only if
    $x_{1}^{(0)}=\cdots=x_{n}^{(0)}=0$.  Thus such a map of rings
    gives a well-defined $m$-jet centered at the origin of $X$ if and
    only if we have
    $$f(t(x_{1}^{(1)}+\ldots+x_{1}^{(m)}t^{m-1}),\ldots,t(x_{n}^{(1)}+
    \ldots+x_{n}^{(m)}t^{m-1})) \in (t^{m+1})$$ for each generator $f$
    of $I$.  Since $I$ is generated by homogeneous degree $d$
    elements, this is equivalent to
    $$f(x_{1}^{(1)}+\ldots+x_{1}^{(m)}t^{m-1},\ldots,x_{n}^{(1)}+
    \ldots+x_{n}^{(m)}t^{m-1}) \in (t^{m+1-d}).$$ But this is the same
    as saying that the ring map
    \begin{align*}
	k[x_{1},\ldots,x_{n}]/I &\rightarrow k[t]/(t^{m-d+1})\\
	x_{i} &\mapsto
	x_{i}^{(1)}+x_{i}^{(2)}t+\ldots+x_{i}^{(m-d+1)}t^{m-d}
    \end{align*}
    is an $(m-d)$-jet of $X$.  Since there are no constraints on the
    variables $x_{i}^{(m-d+2)},\ldots,x_{i}^{(m)}$ for all
    $i=1,\ldots,n$, we see that $\pi_{m}^{-1}(0) \cong
    \mathcal{J}_{m-d}(X) \times \mathbb{A}^{n(d-1)}$.
\end{proof}

\begin{proof}[Proof of Theorem \ref{2x2minors}]
    We will proceed by induction on $m$.  Because we will use
    Proposition \ref{PreimageOfOrigin} to relate $\mathcal{J}_m$ to
    $\mathcal{J}_{m-2}$, we will need base cases for $m=0$ and $m=1$.
    If $m=0$, then $\lfloor \frac{m+1}{2} \rfloor=0$ and
    $\mathcal{J}_{0}(X) \cong X$.  The theorem predicts one component
    of dimension same as that of $X$, which is obvious.  If $m=1$, the
    closed subset $\overline{\pi_{1}^{-1}(X^{reg})}$ is an irreducible
    component of $\mathcal{J}_{1}(X)$ and has dimension twice that of
    $X$.  Note that the singular locus of $X$ is simply the origin, and
    $$k[x_{ij}^{(k)}]/(J_{1}(X)+(x_{ij}^{(0)})) \cong
    k[x_{ij}^{(1)}],$$
    as is easy to see that $J_1 \subseteq (x_{ij}^{(0)}:1 \leq i \leq
    r, 1 \leq j \leq s)$.  So $\pi_{1}^{-1}(0) \cong \spec
    \mathbb{A}^{rs}$.  The condition $r>s \geq 3$ is equivalent to
    \begin{equation}
	rs \geq 2\dim X. \tag{$\star$}
    \end{equation}
    Thus, $\mathcal{J}_{1}(X)$ has two components:
    $\overline{\pi_{1}^{-1}(X^{reg})}$ of dimension $2\dim X$ and
    $\pi_{1}^{-1}(X^{sing})$ of dimension $rs$, as predicted by the
    theorem.  This completes the $m=1$ base case.

    For general $m$, the closed subset
    $\overline{\pi_{m}^{-1}(X^{reg})}$ is irreducible and has
    dimension $(m+1)\dim X$.  On the other hand, Proposition
    \ref{PreimageOfOrigin} tells us $\pi_{m}^{-1}(0) \cong
    \mathcal{J}_{m-2}(X) \times \mathbb{A}^{rs}$.  So by induction,
    $\pi_{m}^{-1}(0)$ has $\lfloor \frac{m-1}{2} \rfloor +1$
    components, where the $q^{th}$ component has dimension
    $$qrs+(m-1-2q)\dim X+rs \hspace{0.5cm}\text{for
    }q=0,\ldots,\lfloor \frac{m-1}{2} \rfloor.$$
    Because $\dim X=r+s-1$, it follows that $\pi_{m}^{-1}(0)$ has
    $\lfloor \frac{m+1}{2} \rfloor$ components, where the $q^{th}$
    component has dimension $$qrs+(m+1-2q)\dim X
    \hspace{0.5cm}\text{where } q=1,\ldots,\lfloor \frac{m+1}{2}
    \rfloor.$$
    Notice that by our assumptions on $r$ and $s$, the minimum value
    is $rs+(m-1)\dim X$, which is greater than or equal to the
    dimension of the component of $\mathcal{J}_{m}(X)$ over the smooth
    part by ($\star$).  Therefore, each of these $\lfloor
    \frac{m+1}{2} \rfloor$ components of $\pi_{m}^{-1}(0)$ is a
    component of $\mathcal{J}_{m}(X)$.  Thus, $\mathcal{J}_{m}(X)$ has
    $\lfloor \frac{m+1}{2} \rfloor +1$ components of dimensions
    $qrs+(m+1-2q)\dim X$ where $q$ ranges from $0$ to $\lfloor
    \frac{m+1}{2} \rfloor$.
\end{proof}

\begin{corollary}\label{LCT}
    With the same assumptions as in Theorem \ref{2x2minors}, the log
    canonical threshold of the pair $(\mathbb{A}^{rs},X)$ is exactly
    $\frac{1}{2}rs$.
\end{corollary}

\begin{proof}
    By applying a result of Musta\c{t}\v{a} \cite[Corollary
    0.2]{MustataSingsPairs}, we obtain
    \begin{align*}
	lct(\mathbb{A}^{rs},X) &= \dim \mathbb{A}^{rs}-\sup_{m \geq 0}
	\frac{\dim \mathcal{J}_{m}(X)}{m+1}\\
	&= rs-\sup_{m \geq 0} \frac{\lfloor \frac{m+1}{2} \rfloor
	rs+(m \text{ mod }2)\dim X}{m+1}\\
	&= rs-\sup_{m \geq 0} \frac{\lfloor \frac{m+1}{2} \rfloor
	rs}{m+1}\\
	&= rs-\frac{1}{2}rs\\
	&= \frac{1}{2}rs.
    \end{align*}
\end{proof}
    
\begin{rmk}
    In her thesis \cite{Johnson}, Johnson produced log canonical 
    thresholds of other determinantal varieties by direct 
    calculations of their log resolutions. 
\end{rmk}    

\begin{rmk}
    The previous three theorems give a host of examples of Gorenstein
    varieties with rational singularities whose jet schemes are not
    irreducible.  In particular, they illustrate that
    Musta\c{t}\v{a}'s result that locally complete intersection
    varieties have rational singularities if and only if their jet
    schemes are irreducible cannot be weakened: we may not replace the
    local complete intersection hypothesis with a Gorenstein
    hypothesis.  For example, taking $r=s=c+3 \geq 4$ in Theorem
    \ref{OddJets} gives a rationally singular Gorenstein variety whose
    odd jet schemes are not irreducible.  Musta\c{t}\v{a} himself also
    gave an example of a toric variety to illustrate this fact
    \cite[Example 4.6]{MustataLCI}.
\end{rmk}

\bibliographystyle{plain}
\bibliography{biblio}   % Use the BibTeX file ``biblio.bib''.

\end{document}